\documentclass[12pt,a4paper]{article}
\usepackage{amsmath,amssymb,amsfonts}
\def\Bbb{\mathbb}

\title{\bf Reducing radicals in the spirit of Euclid}

\author{Kurt Girstmair}
\date{}

\makeatletter
\let\@@maketitle=\maketitle
\def\maketitle{\def\thispagestyle##1{\relax}\@@maketitle}
\makeatother
\newtheorem{theorem}{Theorem}
\newtheorem{prop}{Proposition}

\def\BE{\begin{equation}}
\def\EE{\end{equation}}
\def\BD{\begin{displaymath}}
\def\ED{\end{displaymath}}
\def\BA{\begin{array}}
\def\EA{\end{array}}
\def\BEA{\begin{eqnarray*}}
\def\EEA{\end{eqnarray*}}
\def\BI{\bibitem}

\def\Z{\Bbb Z}
\def\Q{\Bbb Q}

\def\phi{\varphi}

\def\MB{\mbox}
\def\LD{\ldots}
\def\OV{\overline}

\def\DIV{\,|\,}

\def\NI{\noindent}
\def\MN{\medskip\noindent}
\def\STOP{\hfill$\Box$}

\begin{document}
\maketitle

\begin{abstract}
\NI
Let $p$ be an odd natural number $\ge 3$.
Inspired by results from Euclid's {\em Elements}, we express the irrational
\BD
  y=\sqrt[p]{d+\sqrt R},
\ED
whose degree is $2p$, as a polynomial function of irrationals of degrees $\le p$. In certain cases $y$ is expressed by simple radicals. This reduction of the degree
exhibits remarkably regular
patterns of the polynomials involved. The proof is based on hypergeometric summation, in particular, on Zeilberger's
algorithm.
\end{abstract}

\section*{1. Introduction and main result}

In the tenth book of Euclid's elements, Proposition 54, one finds an answer to the following question: Suppose that the biquadratic radical
\BD
\label{1.0}
y=\sqrt{d+\sqrt R}
\ED
is given, where $d$ and $R$ are positive rational numbers, $\sqrt R\not\in\Q$ and $y\not\in\Q(\sqrt R)$.
When is it possible to express $y$ in terms of two simple square roots?
Euclid's answer is as follows: If $d^2-R$ is a square, i.e., $d^2-R=k^2$, $k\in \Q$, then
\BE
\label{1.2}
y=\sqrt{\frac{d+k}2}+\sqrt{\frac{d-k}2}.
\EE
But, to tell the truth, Euclid has no formulas, and his answer is phrased in purely geometric terms. So
(\ref{1.2}) is a modern algebraic interpretation of what Euclid did in the framework of his geometry (see \cite[p. 119]{He}).

Formula (\ref{1.2}) has been considered as an example of {\em denesting} a nested radical (see \cite{Bo}, Th. 1).
Our viewpoint, however, is different, as we illustrate with the example
\BE
\label{1.3}
y=\sqrt[4]{d+\sqrt R}.
\EE
Here we assume that $y$ is an irrational of degree $8$, which means that the polynomial $(Z^4-d)^2-R$ is irreducible.
Provided that $d^2-R=k^4$, $k\in\Q$, we can apply Euclid's result twice and obtain
\BE
\label{1.4}
 y=\sqrt{\sqrt{\frac{d+k^2}8}+\frac k2}+\sqrt{\sqrt{\frac{d+k^2}8}-\frac k2}.
\EE
So (\ref{1.4}) does not denest the nested radical (\ref{1.3}).  Nevertheless, this identity
can be considered as a {\em reduction of the degree}, inasmuch as it expresses an irrational of degree $8$ as a sum of two irrationals of degree $4$.

The present paper is devoted to this kind of reduction in the case of the radical
\BE
\label{1.6}
  y=\sqrt[p]{d+\sqrt R},
\EE
where $p\ge 3$ is an odd natural number. To this end let $K$ be a field of characteristic $0$ and let $d, R\in K\smallsetminus\{0\}$. By $\OV K$ we denote an algebraic closure of $K$. An {\em irrational}
is an element of $\OV K\smallsetminus K$. The {\em degree} of an irrational is the degree of its minimal polynomial over $K$. Let $\sqrt R$ be an irrational. More precisely, we select one of the
two possible values of $\sqrt R$, whereas the other value is denoted by $-\sqrt R$. This meaning of $\sqrt R$ shall be fixed throughout. Suppose that $y$ is an irrational of degree $2p$.
This is the same as saying that the polynomial
\BD
\label{1.8}
  g=(Z^p-d)^2-R
\ED
is irreducible in the polynomial ring $K[Z]$. We are going to reduce $y$ to two irrationals of degree $\le p$ together with $\sqrt R$, i.e., we express $y$ as a
polynomial function (with coefficients in $K$) of these quantities. Surprisingly, this can be done by means of explicit formulas valid for all odd natural numbers $p$.

For this purpose we work with the decomposition
\BE
\label{1.10}
 g=h\cdot h'\enspace \MB{ with }\enspace h=Z^p-d-\sqrt R,\enspace h'=Z^p-d+\sqrt R,
\EE
which takes place in the polynomial ring $K(\sqrt R)[Z]$. Now suppose that $y$ is a zero of $h$, whereas $y'$ $(\in \OV K)$ is a zero of $h'$. We put
\BD
\label{1.12}
 z=yy' \:\MB{ and }\: u=z^{(p-1)/2}(y+y').
\ED
Then $z$ is a $p$th root of
\BE
\label{1.13}
 D=d^2-R,
\EE
i.e., $z^p=D\in K$. By our assumptions, $D\ne 0$.
On the other hand, we will see that $u$ is a zero of the polynomial
\BE
\label{1.14}
  f=D^{(p-1)/2}\sum_{k=0}^{(p-1)/2} \frac{c_{2k+1}}{D^{k}}Z^{2k+1}-2dD^{(p-1)/2}\in K[Z],
\EE
where $D$ is defined by (\ref{1.13}) and $c_{2k+1}$ by
\BE
\label{1.16}
  c_{2k+1}=(-1)^{(p-1)/2-k}\frac p{\frac{p+1}2+k}\binom{\frac{p+1}2+k}{2k+1},
\EE
$k=0,\LD, (p-1)/2$.
Here $\binom{m}{n}$ is the usual binomial coefficient. Note that $f$ is a polynomial of degree $p$ with leading coefficient $c_p=1$.
Define the polynomial $A\in K[Z]$ by
\BE
\label{1.18}
  A=\frac 1{2R}\sum_{k=0}^{(p-1)/2}\frac{a_{2k}}{D^k}Z^{2k}+(-1)^{(p+1)/2}\frac d{2RD} Z
\EE
with
\BE
\label{1.20}
  a_{2k}=(-1)^k\frac{p-1}{\frac{p-1}2+k}\binom{\frac{p-1}2+k}{2k},
\EE
$k=0,\LD, (p-1)/2$.
Then our main result reads as follows.

\begin{theorem}  
\label{t1}
Let $y$, $y'$, $z=yy'$, and $u=z^{(p-1)/2}(y+y')$ be as above. In particular, $\sqrt R\not\in K$, $y, y'$ are irrationals of degree $2p$, and $u$ is a zero of $f$. Then
\BE
\label{1.22}
  \{y,y'\}=\left\{z^{(p+1)/2}\left( \frac u{2D} \pm A(u)\sqrt{R}\right)\right\}.
\EE
\end{theorem} 

\MN
{\em Remarks.} 1. In our setting, $z$ is a $p$th root of $D$ and $u$ a zero of the polynomial $f$ of degree $p$ in $K[Z]$.
Therefore, (\ref{1.22}) reduces the radical
$y$ of (\ref{1.6}), whose degree is $2p$, to $z$, $u$ and $\sqrt{R}$, whose degrees are $\le p$ and $=2$, respectively.

2. Instead of choosing $y$, $y'$ in the above way, one may choose a zero $y$ of $h$ and a $p$th root $z$ of $D$. Then $y'=z/y$ is a zero of the polynomial $h'$.

3. We will see below (Proposition \ref{p1}) that for a given $p$th root $z$ of $D$ and a given zero $u$ of $f$ there is exactly one zero $y$ of h such that for $y'=z/y$
we have $u=z^{(p-1)/2}(y+y')$.

\MN
The proof of Theorem \ref{t1} is based on techniques of hypergeometric summation, in particular, on Zeilberger's algorithm.
We think that one can hardly dispense with these techniques or, in other words, a proof without algorithmic methods of this kind seems to be out of reach.

\MN
{\em Example} 1. We consider the special case $p=5$. Here we have
\BD
   f=Z^5-5DZ^3+5D^2Z-2dD^2 \enspace\MB{ and } A=\frac 1R\left(\frac{Z^4}{2D}-\frac{2Z^2}{D}-\frac{dZ}{2D}+1\right).
\ED
In this case the polynomial $f$ is called {\em DeMoivre's quintic} (see \cite{Brg}, \cite{Sp}).
If, for instance $d=2$ and $R=5$, we obtain $g=Z^{10}-4Z^5-1$, which is irreducible over $K=\Q$. Moreover, $D=-1$, and we can choose $z=-1$ (in the sense of Remark 2). Finally,
$f=Z^5+5Z^3+5Z-4$ is irreducible over $\Q$ and $A=(Z^4+4Z^2+2Z+2)/10$.

\MN
Let $z$ be a $p$th root of $D$ and suppose that $f$ has a zero $u$ in the ground field $K$. By Remark 3, there is a uniquely determined zero $y$ of $h$ such that $u=z^{(p-1)/2}(y+y')$ for $y'=z/y$.
Then (\ref{1.22}) expresses $y$ and $y'$ as polynomial functions of $z$ and $\sqrt R$ with coefficients in $K$. In other words, the nested radicals $y$ and $y'$ are denested in this way.

\MN
{\em Example} 2. Let $K=\Q$ and $p=7$. Put $d=-2158$ and $R=6\cdot 881^2$. Thus, $d+\sqrt R= -2158+881\sqrt 6$ and $D=-2$. In this case the polynomial $g$ equals $Z^{14}+4316Z^7-2$ and is irreducible.
It turns out that $u=4$ is a zero of $f$. We choose $z=-2^{1/7}$, where $2^{1/7}$ is the real $7$th root of $2$. If the zero $y$ of $h$ is determined by $u=z^{(p-1)/2}(y+y')$, $y'=z/y$, then formula (\ref{1.22}) gives
\BD
\label{1.24}
  \{y,y'\}=\left\{2^{4/7}\left(-1\pm \sqrt{6}/2\right)\right\}.
\ED

\MN
The problem of denesting radicals has attracted considerable attention, see, for instance, \cite{Zi}, \cite{Ho}, \cite{Ld1}. Further references can be found in  \cite{Ld}.

Sections 2--4 are devoted to the proof of Theorem \ref{t1}.

Section 5 contains the aforementioned Proposition \ref{p1}. Moreover, we show how to find examples like the above Example 2 in a simple way.

In Section 6 we discuss the question whether $f$ can be replaced by a polynomial of the form $Z^p-a$, $a\in K$, which is the same
as replacing $u$ by the $p$th root of an element of $K$ (Proposition \ref{p2}). In this case our reduction could be replaced by a reduction which also denests the radical $y$.
In addition, we consider two cases where (\ref{1.22}) is equivalent to expressing $y$, $y'$ in terms of a $\Q$-basis of the field $\Q(z,u,\sqrt R)$ (Proposition \ref{p3}). These cases seem to be
generic, inasmuch as most examples fall under them. We also get some insight into the structure of the splitting field of the polynomial $g$.

\section*{2. The polynomial $f$}

We adopt the above notation. In particular, $f$ is defined by (\ref{1.14}) and (\ref{1.16}).
We have to show that $u=z^{(p-1)/2}(y+y')$ is a zero of $f$.

For this purpose we look at the following expansion of the polynomial $X^p+1\in\Q[X]$:
\BE
\label{2.0}
   X^p+1=\sum_{k=0}^{(p-1)/2}C_{p-2k}X^k(X+1)^{p-2k}
\EE
with rational coefficients $C_{p-2k}$. We will see that this expansion exists and is unique. Indeed, if we apply the binomial formula to $(X+1)^{p-2k}$ and consider only the coefficients of
the monomials $X^j$, $j=0,\LD, (p-1)/2$,
we obtain
\BD
  \sum_{j=0}^{(p-1)/2} X^j\sum_{k=0}^j C_{p-2k}\binom{p-2k}{j-k}=1.
\ED
This gives the system of linear equations
\begin{eqnarray}
\label{2.2}
  1&=& C_p, \nonumber\\
  0&=& C_p\binom{p}{1}+C_{p-2}\binom{p-2}{0}, \nonumber\\
  0&=& C_p\binom{p}{2}+C_{p-2}\binom{p-2}{1}+C_{p-4}\binom{p-4}{0},\nonumber\\
\vdots & & \vdots\nonumber\\
  0&=& \sum_{k=0}^j C_{p-2k}\binom{p-2k}{j-k}.
\end{eqnarray}

One immediately sees that the coefficients $C_{p-2k}$, $k=0,\LD, p-1$, are uniquely determined by (\ref{2.2}).

Let $F$ denote the polynomial on the right hand side of (\ref{2.0}). We obtain
\BE
\label{2.4}
 F(y/y')\cdot y'^p=\sum_{k=0}^{(p-1)/2}C_{p-2k}z^k(y+y')^{p-2k}=y^p+y'^p=2d,
\EE
since $y^p=d+\sqrt R$, $y'^p=d-\sqrt R$ (see (\ref{1.10})). If we multiply the identity (\ref{2.4}) by $D^{(p-1)/2}=z^{p(p-1)/2}$, we have
\BD
   2dD^{(p-1/2)}=\sum_{k=0}^{(p-1)/2}C_{p-2k}D^kz^{((p-1)/2)\cdot(p-2k)}(y+y')^{p-2k}.
\ED
In other words, $u=z^{(p-1)/2}(y+y')$ is a zero of $f$, provided that the coefficients $C_{p-2k}$ of (\ref{2.0}) coincide with the
coefficients $c_{p-2k}$ of $f$, $k=0,\LD,(p-1)/2$.

We will show that the coefficients $c_{p-2k}$  satisfy the system (\ref{2.2}) of linear equations.
By (\ref{1.16}),
\BD
 c_{p-2k}=(-1)^k\frac p{p-k}\binom{p-k}{k},
\ED
$k=0,\LD,(p-1)/2$. We recall that $c_p=1$ and observe that the right hand side of (\ref{2.2}) reads, for these values of $c_{p-2k}$,
\BE
\label{2.6}
\sum_{k=0}^j (-1)^k\frac p{p-k}\binom{p-k}{k}\binom{p-2k}{j-k},
\EE
$j=1,\LD,(p-1)/2$.
This is a typical example of a hypergeometric summation (see \cite[chap. 2]{Ko}). The summand of (\ref{2.6}) is defined for all integers $k\in \Z$, since $\binom{p-k}{k}/(p-k)=\binom{p-k-1}{k-1}/k$ for all $k\ne 0$.
In particular, it takes the value $0$ for all $k>j$ and all $k<0$. For our purpose it is advisable to change the summation order, i.e., we consider
\BE
\label{2.8}
\sum_{k=0}^j (-1)^{j-k}\frac p{p-j+k}\binom{p-j+k}{j-k}\binom{p-2j+2k}{k}.
\EE
Let $b_k$ denote the summand of (\ref{2.8}). Hypergeometric summation requires considering $b_{k+1}/b_k$, $k=0,\LD,j$, and $b_0$. We have
\BD
  \frac{b_{k+1}}{b_k}=\frac{(k-j)(k+p-j)}{(k+p-2j+1)(k+1)}.
\ED
Since $b_0=(-1)^j\frac p{p-j}\binom{p-j}{j}$, we obtain that (\ref{2.6}) equals
\BD
  (-1)^j\frac p{p-j}\binom{p-j}{j}\, _2F_1\left(\begin{array}{c|}-j,p-j\\
                                                            p-2j+1
                                              \end{array}\: 1\right),
\ED
$j=1,\LD,(p-1)/2$, where $_2F_1$ denotes Gauss' hypergeometric function. A theorem of Gauss (see \cite[p. 32]{Ko}) says
\BD
   _2F_1\left(\begin{array}{c|}-j,b\\
                              c
               \end{array}\: 1\right)=\frac{(c-b)_j}{(c)_j};
\ED
here $(a)_j$ is given by $(a)_j=a\cdot(a+1)\cdots (a+j-1)$ and $b,c\in\Z$, $c>0$. In our case, $c=p-2j+1>0$ and, thus, $(c)_j\ne 0$. On the other hand, $c-b=-j+1$, which means
$(c-b)_j=0$. Hence the sum (\ref{2.6}) vanishes for all $j=1,\LD,(p-1)/2$.

We remark that the coefficients of the monomials $X^j$, $j=(p+1)/2,\LD,p$, in (\ref{2.0}) also yield the system (\ref{2.2}) for the numbers $C_{p-2k}$.

\section*{3. The basic identity}

We return to Theorem \ref{t1}. Indeed, (\ref{1.22}) is the same as saying
\BE
\label{3.2}
  \pm z^{(p+1)/2}A(u)\sqrt{R}=\frac{y-y'}2.
\EE
For instance, the plus-sign on the left hand side of (\ref{3.2}) gives
\BD
  y=\frac{y+y'}2 +\frac{y-y'}2=\frac u{2z^{(p-1)/2}}+z^{(p+1)/2}A(u)\sqrt{R},
\ED
which obviously has the form of (\ref{1.22}). However, (\ref{3.2}) is equivalent to
\BD
  z^{p+1}A(u)^2R=\frac{(y+y')^2-4yy'}4.
\ED
Since the right hand side of this identity equals $(u^2-4z^p)/(4z^{p-1})$,
we see that (\ref{3.2}) is equivalent to
\BE
\label{3.4}
  4D^2A(u)^2R=u^2-4D.
\EE
We define
\BD
\label{3.6}
 f'=\frac 1{RD^{(p-3)/2}}\sum_{j=0}^{(p-3)/2} \frac{c'_{2j+1}}{D^j}Z^{2j+1}-\frac{2d}{RD^{(p-3)/2}}
\ED
with
\BD
\label{3.8}
  c'_{2j+1}=(-1)^{(p-3)/2-j}\frac {p-2}{\frac{p-1}2+j}\binom{\frac{p-1}2+j}{2j+1}.
\ED
We are going to prove the identity of polynomials
\BE
\label{3.10}
4D^2A^2R=f\cdot f'+Z^2-4D.
\EE
If we insert $u$ for the variable $Z$ in (\ref{3.10}) and observe $f(u)=0$, we obtain (\ref{3.4}).
Hence (\ref{3.10}) can be considered as the fundamental identity of this paper.
The proof of this identity consists in comparing the coefficients of the monomials $Z^m$ on both sides.

The following three cases have to be distinguished. First, $m\in\{0,2\}$, second, $m$ odd, and third, $m=2k, k=2,\LD,p-1$.
The case $m\in\{0,2\}$ may be checked by the reader. In the remaining cases, we write $\alpha_m$ for the coefficient of $Z^m$ on the left hand side of
(\ref{3.10}) and $\beta_m$ for the coefficient of $Z^m$ on the right hand side. In view of (\ref{1.18}), (\ref{1.20}),
we obtain
\BD
 RD^{k-1}\alpha_{2k+1}=(-1)^{\frac{p+1}2+k}2d\cdot\frac{p-1}{\frac{p-1}2+k}\binom{\frac{p-1}2+k}{2k}
\ED
and
\BD
 RD^{k-1}\beta_{2k+1}=(-1)^{\frac{p+1}2+k}2d\cdot\left(\frac p{\frac{p+1}2+k}\binom{\frac{p+1}2+k}{2k+1}
             -\frac {p-2}{\frac{p-1}2+k}\binom{\frac{p-1}2+k}{2k+1}\right),
\ED
$k=0,\LD,(p-1)/2$ (observe that $\binom{p-1}{p}=0$).
Using elementary identities of binomial coefficients, one sees that $\alpha_{2k+1}=\beta_{2k+1}$.

The remaining case is the most difficult one. We obtain
\BE
\label{3.12}
   RD^{k-2}\alpha_{2k}=(-1)^k\sum_{j=0}^k\frac{(p-1)^2}{(\frac{p-1}2+j)(\frac{p-1}2+k-j)}\binom{\frac{p-1}2+j}{2j}\binom{\frac{p-1}2+k-j}{2(k-j)}
\EE
and
\BE
\label{3.14}
 RD^{k-2}\beta_{2k}=(-1)^k\sum_{j=0}^{k-1}\frac{p(p-2)}{(\frac{p+1}2+j)(\frac{p-1}2+k-j-1)}\binom{\frac{p+1}2+j}{2j+1}\binom{\frac{p-1}2+k-j-1}{2(k-j-1)+1}
\EE
for $k=2,\LD, (p-1)$. Here we observe that the left binomial coefficient in (\ref{3.12}) as well as in (\ref{3.14}) equals $0$ if $j>(p-1)/2$. In the same way
the right binomial coefficient vanishes in both identities if $k-j>(p-1)/2$. We also observe that the sum on the right hand side of (\ref{3.14}) may be extended to the upper bound $k$ (instead of $k-1$)
since the respective summand is $0$.
In the next section we show that $\alpha_{2k}=\beta_{2k}$, $k=2,\LD, p-1$.

\section*{4. Two hypergeometric summations}

In this section we denote the right hand side of (\ref{3.12}) by $s_k$ and the right hand side of (\ref{3.14}) by $t_k$ for $k=0,\LD, p-1$. Further, we introduce
\BD
\label{4.2}
    u_k=\frac{(-1)^k (p-1)}{k}\binom{p+k-2}{2k-1},
\ED
$k=1,\LD,p-1$. We will show that $s_k=u_k$ and $t_k=u_k$ for all $k=2,\LD,p-1$. In this way we also exhibit the value of the coefficient $\alpha_{2k}$ of $4D^2A^2R$.

Zeilberger's algorithm yields recursion formulas for $s_k$ and $t_k$ (see \cite[chap. 7]{Ko}). In the case of $s_k$, this formula reads
\BE
\label{4.4}
  (4k^2+6k+2)s_{k+1}+(-k^2-2p+1+p^2)s_k=0.
\EE
Note that $4k^2+6k+2\ne 0$ for $k\ge 0$.
Now $u_k\ne 0$ for $k=1,\LD,p-2$ and
\BD
 \frac{u_{k+1}}{u_k}=-\frac{-k^2-2p+1+p^2}{4k^2+6k+2}.
\ED
In particular, $u_k$ satisfies formula (\ref{4.4}). Moreover, $u_1=s_1=-(p-1)^2$. Hence $u_k=s_k$ for all $k=1,\LD,p-1$.

In the case of $t_k$ Zeilberger's algorithm yields
\BD
\label{4.6}
  a\cdot t_{k+2}+b\cdot t_{k+1} +c\cdot t_k=0,
\ED
with
\BEA
  a&=&16k^3+64k^2+76k+24,\: b\enspace =\enspace -8k^3-12k^2-8pk+4p^2k+2p^2-4p+2,\\
  c&=&k^3-k^2-p^2k+2pk-k+p^2-2p+1.
\EEA
Again, $a\ne 0$ for all $k\ge 0$.
It is not hard to check that
\BD
 a\cdot\frac{u_{k+2}}{u_k}+b\cdot\frac{u_{k+1}}{u_k}+c=0
\ED
for $k=2,\LD, p-3$. Since $t_2=u_2=p(p-1)^2(p-2)/12$ and $t_3=u_3=-p(p-1)^2(p-2)(p-3)(p+1)/360$, we obtain
$u_k=t_k$ for all $k=2,\LD,p-1$.

\section*{5. Some additional observations}

Let the above assumptions hold, in particular, $d\ne 0$ and $\sqrt R\not\in K$.

\begin{prop} 
\label{p1}

Let $z$ be an arbitrary $p$th root of $D$ and $u$ a zero of $f$. Then there is a uniquely determined zero $y$ of $h$ such that
for $y'=z/y$ we have $u=z^{(p-1)/2}(y+y')$.

\end{prop} 

\MN
{\em Proof.}
Let $V(h)$ and $V(f)$ denote the sets of the zeros (in $\OV K$) of $h$ and $f$, respectively. Then the map
\BE
\label{5.2}
V(h)\to V(f):y\mapsto z^{(p-1)/2}(y+z/y)
\EE
is well defined (recall Remark 2 in Section 1). This map is injective. Indeed, if $u\in V(f)$ equals $z^{(p-1)/2}(y+z/y)$ for some $y\in V(h)$, then the quadratic equation
\BE
\label{5.4}
z^{(p-1)/2}(x+z/x)=u
\EE
has at most two solutions, namely, $x=y$ and $x=z/y$. However, $z/y$ is a zero of $h'$, and $h$ and $h'$ have no common zero since $\sqrt R\ne 0$. Hence $y$ is the only solution of this equation in $V(h)$.

Since $h$ has no multiple zeros, the set $V(h)$ has $p$ elements and $V(f)$ at most $p$. By the injectivity, $V(f)$ also has $p$ elements and the map of (\ref{5.2}) is bijective. This proves our assertion.
\STOP

Our next aim is a simple construction of examples like Example 2 in Section 1. Suppose that $D\in K\smallsetminus\{0\}$ and $u\in K$ are given. Then the equation $f(u)=0$ holds if, and only if,
\BD
d=\frac 12\sum_{j=0}^{(p-1)/2}\frac{c_{2j+1}}{D^j}u^{2j+1}
\ED
(recall (\ref{1.14}), (\ref{1.16})). For instance, if $D=-2$ and $u=4$, we obtain $d=-2158$ (see the aforementioned example). Then we determine $R$ by $R=d^2-D$. In the case $K=\Q$ it frequently happens that $\sqrt R\not\in\Q$ and
that $g=(Z^p-d)^2-R$ is irreducible in $\Q[Z]$. Choose a $p$th root $z$ of $D$ (in $\OV K)$). Then the zero $y$ of $h$ is uniquely determined as a solution of the quadratic equation (\ref{5.4});
and $y$ and $y'=z/y$ can be obtained by (\ref{1.22}). In general, however, it is simpler to obtain $y,y'$ as solutions of the quadratic equation --- provided that $u$ is known.  In this case
\BE
\label{5.6}
  \{y,y'\}=\left\{\frac 1{2z^{(p-1)/2}}\left(u \pm \sqrt{u^2-4z^p}\right)\right\}.
\EE
Since $u$ and $z^p=D\in K$, we have $u^2-4z^p\in K$. Of course, (\ref{5.6}) is equivalent to (\ref{1.22}) in this context.

\section*{6. Two further results}

The question arises whether the zero $u$ of $f$ can be replaced by the $p$th root of an element of $K$. This would imply that the splitting field $L$ of $f$ is also the splitting field of a polynomial $P$ of the form
$P=Z^p-a$ for some $a\in K$.

In order to obtain a partial answer, we suppose that $p\ge 3$ is a prime and $K=\Q$. Further, we assume that $f$ is irreducible over $\Q$.
If $L$ is also the splitting field of $P$, then $P$ is irreducible over $\Q$ and $L$ contains a primitive $p$th root of unity $\zeta_p$. Indeed, $L=\Q(v,\zeta_p)$, where $v$ is a zero of $P$.

\begin{prop} 
\label{p2}

As above, let $K=\Q$, $p\ge 3$ a prime and $f$ irreducible in $\Q[Z]$. If $\Q(\sqrt R)\ne \Q(\sqrt{(-1)^{(p-1)/2}p})$, then $\zeta_p\not\in L$. In particular, $L$ is not the splitting field
of a polynomial $P$ as above.

\end{prop} 

\MN
{\em Proof.}
Let $y$ be a zero of $h$, $y'$ a zero of $h'$, and $z=yy'$. The bijection (\ref{5.2}) implies that the numbers $u_k=z^{(p-1)/2}(y\zeta_p^k+y'\zeta_p^{-k})$, $k=0,\LD,p-1$, are exactly the zeros of $f$. The Lagrange resolvent
\BD
  \sum_{k=0}^{p-1}\zeta_p^{-k}u_k=z^{(p-1)/2}yp+ z^{(p-1)/2}y'\sum_{k=0}^{p-1}\zeta_p^{-2k}=z^{(p-1)/2}yp
\ED
shows $z^{(p-1)/2}y\in L(\zeta_p)$ and, thus, $D^{(p-1)/2}y^p=D^{(p-1)/2}(d+\sqrt R)\in L(\zeta_p)$.
In particular, $\sqrt R\in L(\zeta_p)$. Observe that all elements of $L$ of a degree different from $p$ are contained in a uniquely determined subfield $L_1$, whose degree (over $\Q)$ divides $p-1$ (see \cite[p. 163]{Hu}).
Suppose that $\zeta_p\in L$. Then $\zeta_p\in L_1$, and, therefore, $L_1=\Q(\zeta_p)$ since $\zeta_p$ has the degree $p-1$ over $\Q$. Moreover, $\sqrt R\in L_1$.
The field $\Q(\zeta_p)$ has a cyclic Galois group of order $p-1$ over $\Q$. Accordingly, it contains a uniquely determined
quadratic subfield, namely, $\Q(\sqrt{(-1)^{(p-1)/2}p})$ (see \cite[p. 71]{Ir}). This implies $\Q(\sqrt R)=\Q(\sqrt{(-1)^{(p-1)/2}p})$.
\STOP

\MN
Next we investigate the connection between formula (\ref{1.22}) and $\Q$-bases.

 \begin{prop} 
\label{p3}

As above, let $K=\Q$, $p\ge 3$ a prime and $f$ irreducible in $\Q[Z]$. In addition, suppose that $\zeta_p$ is not contained in the splitting field of $f$.

{\rm (a)} If $D=z^p$ for some $z\in \Q$, then $u^k\sqrt R^{\,l}$, $k=0,\LD,p-1$, $l=0,1$, is a $\Q$-basis of $\Q(y)=\Q(y')=\Q(u,\sqrt R)$.

{\rm (b)} Suppose that $D$ does not have this form. Let $z$ be a $p$th root of $D$. Then $z^ju^k\sqrt R^{\,l}$, $j,k=0,\LD,p-1$, $l=0,1$, is a $\Q$-basis of $\Q(z,u,\sqrt R)$.

In both cases, {\rm (\ref{1.22})} expresses $y,y'$ in terms of the respective basis.
\end{prop} 

\MN
{\em Proof.}
Assertion (a) is obvious.
Suppose that $D$ does not have the form of (a). Then $Z^p-D$ is irreducible
over $\Q$ (see \cite[p. 221]{La}). Let $z$ be a $p$th root of $D$. Suppose that $Z^p-D$ has a zero in $\Q(u)$. This means that $z\zeta_p^k\in\Q(u)$ for some $k\in\{0,\LD,p-1\}$. But then $\Q(u)=\Q(z\zeta_p^k)$,
because both $f$ and $Z^p-D$ are irreducible. Since the splitting field $L$ of $f$ is a Galois extension of $\Q$, this implies that $Z^p-D$ splits into linear factors over $L$. In particular, $\zeta_p\in L$, which we have excluded.
Hence $Z^p-D$ has no zero in $\Q(u)$, and, by the cited argument, $Z^p-D$ is irreducible over $\Q(u)$. Accordingly, $\Q(z,u)$ has the degree $p^2$ over $\Q$. This field does not contain $\sqrt R$ for reasons
of degree. Altogether, the field $\Q(z,u,\sqrt R)$ has the degree $2p^2$ over $\Q$  and the $\Q$-basis of (b). Now it is clear that (\ref{1.22}) expresses $y, y'$ in terms of this basis.
\STOP

\MN
Note that in case (b) neither $y\in\Q(u,\sqrt R)$ nor $y\in\Q(z,\sqrt R)$. Observe, further, that in this case only $(p+1)/2+2$ of a total of $2p^2$ basis vectors actually occur in (\ref{1.22}).
It seems that Propositions \ref{p2} and \ref{p3} cover the generic case, i.e., most examples satisfy the assumptions of these propositions.

Let us briefly look at the splitting field $M$ of the polynomial $g$ in case (b). If $L$ denotes the splitting field of $f$, then $M$ is the composite of the Galois extensions $\Q(z,\zeta_p)$ and $L$ of $\Q$
(recall that $\sqrt R\in L(\zeta_p)$). The degrees of these extensions over $\Q$ are $p(p-1)$ and $pq$, $q\DIV p-1$, respectively.
The structure of their Galois groups is well known (see \cite[p. 163]{Hu}). The intersection of these Galois extensions is a subfield of $\Q(\zeta_p)$.


\vspace{0.5cm}
\noindent
Kurt Girstmair            \\
Institut f\"ur Mathematik \\
Universit\"at Innsbruck   \\
Technikerstr. 13/7        \\
A-6020 Innsbruck, Austria \\
Kurt.Girstmair@uibk.ac.at

\end{document}